\documentclass[12pt]{amsart}
\usepackage{amsmath,amsfonts,amssymb,amscd}
\usepackage{latexsym}

\font\dfont=cmbx10 at 11pt   %%% FOR DEFIN

\title [Perfect representations]
{Irreducibility of perfect representations of
double affine Hecke algebras}
\author[Ivan Cherednik]{Ivan Cherednik $^\dag$} 
\date{June 6, 2004}

\thanks{$^\dag$ \ Partially supported by NSF grant 
DMS-0200276}

\address[I. Cherednik]{Department of Mathematics, UNC 
Chapel Hill, North Carolina 27599, USA\\
chered@email.unc.edu}

 \def\bysame{{\bf --- }}  
 \def\~{{\bf --}}

\newcommand{\Z}{{\mathbb Z}}
\newcommand{\Q}{{\mathbb Q}}
\newcommand{\N}{{\mathbb N}}

\newcommand{\R}{{\mathbb R}}

\def\HH{\mbox{${\mathcal H}$\kern-5.2pt${\mathcal H}$}}

\newtheorem{theorem}{Theorem}[section]
\newtheorem{proposition}[theorem]{Proposition}
\newtheorem{definition}[theorem]{Definition}
\newtheorem{lemma}[theorem]{Lemma}

\newtheorem{theorem }{Theorem}[section]
\newtheorem{proposition }[theorem]{Proposition}
\newtheorem{definition }[theorem]{Definition}
\newtheorem{lemma }[theorem]{Lemma}
\newtheorem{corollary }[theorem]{Corollary}
\newtheorem{notation }[theorem]{Notation}
\newtheorem{remark }[theorem]{Remark}
\newtheorem{example }[theorem]{Example}

\newtheorem{ theorem}{Theorem}[section]
\newtheorem{ proposition}[theorem]{Proposition}
\newtheorem{ definition}[theorem]{Definition}
\newtheorem{ lemma}[theorem]{Lemma}
\newtheorem{ corollary}[theorem]{Corollary}
\newtheorem{ notation}[theorem]{Notation}
\newtheorem{ remark}[theorem]{Remark}
\newtheorem{ example}[theorem]{Example}

\def\for{\  \hbox{ for } \ }
\def\iif{ \ \hbox{ if } \ }

\def\where{\  \hbox{ where } \ }
\def\and{\  \hbox{ and } \ }
\def\and{\  \hbox{ and } \ }

\def\equal{\stackrel{\,\mathbf{def}}{= \kern-3pt =}}

\def\la{\lambda}

\def\om{\omega}

\def\al{\alpha}

\def\ga{\gamma}
\def\ep{\epsilon}

\def\de{\delta}

\def\si{\sigma}

\def\Ga{\Gamma}
\def\ze{\zeta}

\def\vth{{\vartheta}}

\def\tal{\tilde{\alpha}}

\def\tGa{\tilde{\Gamma}}

\def\tw{\widetilde w}
\def\tW{\widetilde W}

\def\tz{\tilde z}
\def\tb{\tilde b}

\def\tR{\tilde R}

\def\hw{\widehat{w}}
\def\hW{\widehat{W}}

\def\hv{\hat{v}}

\def\0{\mathbf{0}}

\def\çF{\mathcal{F}}

\def\h{\mathcal{H}}

\def\e{\mathcal{E}}
\def\v{\mathcal{V}}

\def\lan{\langle}
\def\llb{(\!(}
\def\ran{\rangle}
\def\rrb{)\!)}
\def\lng{\hbox{\tiny lng}}
\def\sht{\hbox{\tiny sht}}

\newcommand{\sq}{\phantom{1}\hfill$\qed$}

\newcommand{\sgn}{\mbox{sgn}}

\def\HH{\mathfrak{H}}

\def\HH{\hbox{${\mathcal H}$\kern-5.2pt${\mathcal H}$}}

\font\smm=msbm10 at 12pt 
\def\symbol#1{\hbox{\smm #1}}
\def\lsmash{{\symbol n}}

\def\#{\sharp}

\begin{document}
\maketitle
{\small
\tableofcontents
}

%$\mathfrak{a,b,c,d,e,f,g,h,i,j,k,l,m,n,o,p,q,r,s,t,e^u,e^v,e^w}$
%\vfil

In the paper, we prove that
the quotient of the polynomial representation
of the double affine Hecke algebra (DAHA) by the
radical of the duality pairing is always irreducible
(apart from the roots of unity) provided that it is 
finite dimensional. We also find necessary and
sufficient conditions for the radical to be zero, 
which is a $q$-generalization of Opdam's formula for the
{\em singular} $k$-parameters with the multiple 
zero-eigenvalue of the corresponding Dunkl operators. 

Concerning the terminology,
{\em perfect modules} in the paper are finite
dimensional possessing a non-degenerate {\em duality pairing}. 
The latter induces the canonical 
{\em duality anti-involution} of DAHA.
Actually, it suffices to assume
that the pairing is {\em perfect}, i.e. identifies the
module with its dual as a vector space, but
we will stick to the finite dimensional case. 

We also assume that perfect modules are 
{\em spherical}, i.e., quotients of
the polynomial representation of DAHA,
and invariant under the {\em projective action of 
$PSL(2,\Z).$} We do not impose
the semisimplicity in contrast to  
\cite{C12}. The irreducibility theorem  
from this paper is stronger and
at the same time the proof is simpler than 
those in \cite{C12}. 

The irreducibility follows from the projectiv
$PSL(2,\Z)$-action which readily results
from the $\tau_-$-invariance.  
The latter always holds if $q$ is not a root of unity.
At roots of unity, it is true for special $k$ only.
We do not give
in the paper necessary and sufficient conditions for 
the $\tau_-$-invariance as $q$ is a root of unity. 
Generally, it is not difficult to check 
(if it is true).

\smallskip

The polynomial representation has the canonical duality paring. 
It is defined in terms of the difference-trigonometric
Dunkl operators, similar
to the rational case where the differential-rational operators
are used, and involves 
the evaluation at $q^{-\rho_k}$ instead of the value at zero.
The quotient of the
polynomial representation
by the radical $Rad$ of this pairing is a universal 
{\em quasi-perfect} representation. By the latter, we mean 
a DAHA-module with a
non-degenerate but maybe non-perfect 
duality paring.

The polynomial representation, denoted by $\v$ in the
paper, is quasi-perfect and irreducible
for generic values of the DAHA-parameters $q,t.$
It is also $Y$-semisimple, i.e.,
there exists a basis of eigenvectors
of the $Y$-operators, and has the simple $Y$-spectrum 
for generic $q,t.$ 
  
The radical $Rad$ is nonzero when $q$ is a root
of unity or as $t=\varsigma\, q^k$ for {\em special fractional\ } $k$
and proper roots of unity $\varsigma.$

We give an example of {\em reducible} $\v$ which has no
radical ($B_n$). The complete list will be presented in
the next paper.

\medskip
{\bf Semisimplicity.}
Typical examples of $Y$-semisimple perfect
representations are the {\em nonsymmetric Verlinde algebras},
generalizing the Verlinde algebras. The latter describe the fusion
of the integrable representations of the Kac-Moody algebras,
and, equivalently, the reduced category of representations
of quantum groups at roots of unity. The third interpretation
is via factors/subfactors. Generally, these algebras appear
in terms of the vertex operators (coinvariants) associated
with Kac-Moody or Virasoro-type algebras.

There are at least two important reasons to drop the 
semisimplicity constraint:

\smallskip
First, it was found recently that the fusion procedure
for a certain Virasoro-type algebra leads to a non-semisimple 
variant of the Verlinde algebra. 
As a matter of fact, there are no general 
reasons to expect semisimplicity in the 
{\em massless conformal field theory.} The positive
definite inner product in the Verlinde algebra, which
guarantees the semisimplicity, is
given in terms of the masses of the points/particles.   

Second, non-semisimple representations of DAHA are
expected to appear when the whole category of 
representations of {\em Lusztig's
quantum group} at roots of unity is considered. Generally,
non-spherical representations could be necessary. However 
the anti-spherical (Steinberg-type) representations, which are 
spherical constructed for $t^{-1}$ in place of $t,$ 
are expected to play an important role. 
\smallskip

The simplest {\em non-semisimple} example
at roots of unity ($A_1$) is considered at the end of the paper  
in detail.
\smallskip

Concerning the necessary and sufficient condition
for the radical of $\v$ to be
nonzero, it readily follows from the evaluation formula 
for the nonsymmetric Macdonald polynomials \cite{C4}. 
This approach does require the $q,t$-setting because
the evaluation formula collapses in the limit.
Cf. \cite{DO}, Section 3.2.

The method from \cite{O6} (see also \cite{DJO} and
\cite{Je1}) based on the shift operator is also possible,
and even becomes simpler with $q,t$ than in the 
rational/trigonometric case.
It will be demonstrated in the next paper.
The definition of the radical of the polynomial representation
is due to Opdam in the rational case.
See, e.g., \cite{DO}.  
In the $q,t$-case, the radical was introduced in \cite{C3,C4}.   

\medskip
{\bf Rational limit.}
Interestingly, the quotient of $\v$ by the radical is always 
irreducible for the {\em rational DAHA.}
The justification is immediate and goes as follows.

This quotient has the 
zero-eigenvalue (no other eigenvalues appear
in the rational setting) of multiplicity one. Any its 
proper submodule will generate at least one additional 
zero-eigenvector, which is impossible. 
\smallskip

The DAHA and its rational degeneration are 
connected by exp - log maps of some kind \cite{C29}, 
but these maps are of analytic nature in the infinite 
dimensional case
and cannot be directly applied to the polynomial representation. 

Generally, the $q,t$-methods are simpler in
many aspects than those in the rational degeneration 
thanks to the existence of the Macdonald polynomials and their
analytic counterparts. It is somewhat 
similar to the usage of the unitary invariant scalar product
in the theory of {\em compact} Lie groups
vs. the abstract theory of Lie algebras.
The $q,t$-generalization of Opdam's formula for singular $k$ and 
the theory of perfect representations
are typical examples in favor of the $q,t$-setting. 
However, with the irreducibility of the universal quasi-perfect quotient 
of the polynomial representation, it is the other way round. 
\smallskip

My guess is that it happens because the $q,t$-polynomial 
representation contains more information than  
could be seen after the rational degeneration. 
I mean mainly the semisimplicity which
do not exist in the rational theory and can be incorporated only
if the rational DAHA is extended by the "first jet" towards
$q$ (not published). 
\smallskip

It must be mentioned here that the rational theory is for
complex reflection groups. The $q,t$-theory is mainly
about the crystallographic groups. Not all complex reflection
groups have affine and double affine extensions.   
\smallskip

I thank A.~Garsia, E.~Opdam, and N.~Wallach for useful
discussions. I would like to thank
UC at San Diego and IML (Luminy) for
the kind invitations.

\medskip
%\vfill
%\vskip 0.2cm
\setcounter{equation}{0}
\section{Affine Weyl groups}
\setcounter{equation}{0}

Let $R=\{\al\}   \subset \R^n$ be a root system of type 
$A,B,...,F,G$
with respect to a euclidean form $(z,z')$ on $\R^n 
\ni z,z'$,
$W$ the {\dfont Weyl group} 
generated by the reflections $s_\al$,
$R_{+}$ the set of positive  roots ($R_-=-R_+$), 
corresponding to (fixed) simple roots  
roots $\al_1,...,\al_n,$ 
$\Ga$ the Dynkin diagram  
with $\{\al_i, 1 \le i \le n\}$ as the vertices. 

We will also use 
the dual roots (coroots) and the dual root system: 
$$R^\vee=\{\al^\vee =2\al/(\al,\al)\}.$$

The root lattice and the weight lattice are: 
\begin{align}
& Q=\oplus^n_{i=1}\Z \al_i \subset P=\oplus^n_{i=1}\Z \om_i, 
\notag 
\end{align}
where $\{\om_i\}$ are fundamental weights:
$ (\om_i,\al_j^\vee)=\de_{ij}$ for the 
simple coroots $\al_i^\vee.$

Replacing $\Z$ by $\Z_{\pm}=\{m\in\Z, \pm m\ge 0\}$ we obtain
$Q_\pm, P_\pm.$
Note that $Q\cap P_+\subset Q_+.$ Moreover, each $\om_j$ has all 
nonzero coefficients (sometimes rational) when expressed 
in terms of 
$\{\al_i\}.$
Here and further see  \cite{Bo}.  

The form will be normalized
by the condition  $(\al,\al)=2$ for the 
{\em short} roots. Thus,

\centerline{ 
$\nu_\al\equal (\al,\al)/2$ is either $1,$ or $\{1,2\},$ or 
$\{1,3\}.$ }
 
We will use the 
notation $\nu_{\lng}$ for the long roots ($\nu_{\sht}=1$).

\smallskip
Let  $\vth\in R^\vee $ be the {\dfont maximal positive 
coroot}. Considered as a root 
(it belongs to $R$ because of the choice of normalization)
it is maximal among all short positive roots of $R.$

Setting
$\nu_i\ =\ \nu_{\al_i}, \ 
\nu_R\ = \{\nu_{\al}, \al\in R\},$ one has 
\begin{align}\label{partialrho}
&\rho_\nu\equal (1/2)\sum_{\nu_{\al}=\nu} \al \ =
\ \sum_{\nu_i=\nu}  \om_i, \hbox{\ where\ } \al\in R_+,
\ \nu\in\nu_R.
\end{align}
Note that $(\rho_\nu,\al_i^\vee)=1$ as $\nu_i=\nu.$
We will call $\rho_\nu$ {\dfont partial} $\rho.$
\smallskip

{\bf Affine roots.}
The vectors $\ \tal=[\al,\nu_\al j] \in 
\R^n\times \R \subset \R^{n+1}$ 
for $\al \in R, j \in \Z $ form the 
{\dfont affine root system} 
$\tR \supset R$ ($z\in \R^n$ are identified with $ [z,0]$).  
We add $\al_0 \equal [-\vth,1]$ to the simple
 roots for the 
maximal short root $\vth$.
The corresponding set $\tR$ of positive roots coincides with
$R_+\cup \{[\al,\nu_\al j],\ \al\in R, \ j > 0\}$.

We complete the Dynkin diagram $\Ga$ of $R$  
by $\al_0$ (by $-\vth$ to be more
exact). The notation is $\tGa$. One can obtain it from the
completed Dynkin diagram for $R^\vee$ from \cite{Bo} 
reversing the arrows. The number of laces between $\al_i$ and
$\al_j$ in $\tGa$ is denoted by $m_{ij}.$

The set of
the indices of the images of $\al_0$ by all 
the automorphisms of $\tGa$ will be denoted by $O$ 
($O=\{0\} \for E_8,F_4,G_2$). Let $O'={r\in O, r\neq 0}$.
The elements $\om_r$ for $r\in O'$ are the so-called minuscule
weights: $(\om_r,\al^\vee)\le 1$ for
$\al \in R_+$.

Given $\tal=[\al,\nu_\al j]\in \tR,  \ b \in B$, let  
\begin{align}
&s_{\tal}(\tz)\ =\  \tz-(z,\al^\vee)\tal,\ 
\ b'(\tz)\ =\ [z,\ze-(z,b)]
\label{ondon}
\end{align}
for $\tz=[z,\ze] \in \R^{n+1}$.

The {\dfont affine Weyl group} $\tW$ is generated by all $s_{\tal}$
(we write $\tW = \lan s_{\tal}, \tal\in \tR_+\ran)$. One can take
the simple reflections $s_i=s_{\al_i}\ (0 \le i \le n)$ 
as its generators and introduce the corresponding notion of the  
length. This group is the semidirect product $W\lsmash Q'$ of 
its subgroups $W=$ $\lan s_\al,
\al \in R_+\ran$ and $Q'=\{a', a\in Q\}$, where
\begin{align}
& \al'=\ s_{\al}s_{[\al,\nu_{\al}]}=\ 
s_{[-\al,\nu_\al]}s_{\al}\for 
\al\in R.
\label{ondtwo}
\end{align}

The 
{\dfont extended Weyl group} $ \hW$ generated by $W\and P'$
(instead of $Q'$) is isomorphic to $W\lsmash P'$:
\begin{align}
&(wb')([z,\ze])\ =\ [w(z),\ze-(z,b)] \for w\in W, b\in P.
\label{ondthr}
\end{align}
From now on,  $b$ and $b',$ $P$ and $P'$ will be identified.

Given $b\in P_+$, let $w^b_0$ be the longest element
in the subgroup $W_0^{b}\subset W$ of the elements
preserving $b$. This subgroup is generated by simple 
reflections. We set
\begin{align}
&u_{b} = w_0w^b_0  \in  W,\ \pi_{b} =
b( u_{b})^{-1}
\ \in \ \hW, \  u_i= u_{\om_i},\pi_i=\pi_{\om_i},
\label{xwo}
\end{align}
where $w_0$ is the longest element in $W,$
$1\le i\le n.$

The elements $\pi_r\equal\pi_{\om_r}, r \in O'$ and
$\pi_0=\hbox{id}$ leave $\tGa$ invariant 
and form a group denoted by $\Pi$, 
 which is isomorphic to $P/Q$ by the natural 
projection $\{\om_r \mapsto \pi_r\}$. As to $\{ u_r\}$,
they preserve the set $\{-\vth,\al_i, i>0\}$.
The relations $\pi_r(\al_0)= \al_r= ( u_r)^{-1}(-\vth)$ 
distinguish the
indices $r \in O'$. Moreover (see e.g., \cite{C12}):
\begin{align}
& \hW  = \Pi \lsmash \tW, \where
  \pi_rs_i\pi_r^{-1}  =  s_j \iif \pi_r(\al_i)=\al_j,\ 
 0\le j\le n.
\end{align}

Setting
$\hw = \pi_r\tw \in \hW,\ \pi_r\in \Pi, \tw\in \tW,$
the length $l(\hw)$ 
is by definition the length of the reduced decomposition 
$\tw= $ $s_{i_l}...s_{i_2} s_{i_1} $
in terms of the simple reflections 
$s_i, 0\le i\le n.$ 

The length can be 
also defined as the 
cardinality $|\la(\hw)|$
of  
$$
\la(\hw)\equal\tR_+\cap \hw^{-1}(\tR_-)=\{\tal\in \tR_+,\ 
\hw(\tal)\in \tR_-\},\ 
\hw\in \hW.
$$ 

\smallskip
{\bf Reduction modulo $W$.} The following proposition 
is from \cite{C4}.
It generalizes the construction of the elements 
$\pi_{b}$ for $b\in P_+.$ 

\begin{proposition} \label{PIOM}
 Given $ b\in P$, there exists a unique decomposition 
$b= \pi_b  u_b,$
$ u_b \in W$ satisfying one of the following equivalent conditions:

{i) \ \ } $l(\pi_b)+l( u_b)\ =\ l(b)$ and 
$l( u_b)$ is the greatest possible,

{ii)\  }
$ \la(\pi_b)\cap R\ =\ \emptyset$. 
Moreover, $ u_b(b)
\equal b_-\in P_-=-P_+$
is a unique element
from $P_{-}$ which belongs to the orbit $W(b).$ \sq 
\end{proposition}

For $\tal=[\al,\nu_\al j]\in \tR_+,$ one has:
\begin{align}
\la(b) = \{ \tal,\  &( b, \al^\vee )>j\ge 0 \iif \al\in R_+,
\label{xlambi}\\ 
&( b, \al^\vee )\ge j> 0 \iif \al\in R_-\},
\notag \\  
\la(\pi_b) = \{ \tal,\ \al\in R_-,\ 
&( b_-, \al^\vee )>j> 0 
\iif  u_b^{-1}(\al)\in R_+,
\label{lambpi} \\    
&( b_-, \al^\vee )\ge j > 0 \iif   
u_b^{-1}(\al)\in R_- \}, \notag \\ 
\la(u_b) = \{ \al\in R_+,\  &(b,\al^\vee)> 0 \}.
\label{laumin}  
\end{align}

\medskip
\section{Double Hecke algebras}
\setcounter{equation}{0}
By  $m,$ we denote the least natural number 
such that  $(P,P)=(1/m)\Z.$  Thus
$m=2 \for D_{2k},\ m=1 \for B_{2k} \and C_{k},$
otherwise $m=|\Pi|$.

The double affine Hecke algebra depends 
on the parameters 
$q, t_\nu,\, \nu\in \{\nu_\al\}.$ The definition ring is 
$\Q_{q,t}\equal$
$\Q[q^{\pm 1/m},t^{\pm 1/2}]$ formed by the
polynomials in terms of $q^{\pm 1/m}$ and  
$\{t_\nu^{\pm 1/2} \}.$
We set
\begin{align}
&   t_{\tal} = t_{\al}=t_{\nu_\al},\ t_i = t_{\al_i},\ 
q_{\tal}=q^{\nu_\al},\ q_i=q^{\nu_{\al_i}},\notag\\ 
&\where \tal=[\al,\nu_\al j] \in \tR,\ 0\le i\le n.
\label{taljx}
\end{align}

It will be convenient to use the parameters
$\{k_\nu\}$ together with  $\{t_\nu \},$ setting
$$
t_\al=t_\nu=q_\al^{k_\nu} \for \nu=\nu_\al, \and
\rho_k=(1/2)\sum_{\al>0} k_\al \al.
$$

For pairwise commutative $X_1,\ldots,X_n,$    
\begin{align}
& X_{\tb}\ =\ \prod_{i=1}^nX_i^{l_i} q^{ j} 
\iif \tb=[b,j],\ \hw(X_{\tb})\ =\ X_{\hw(\tb)}.
\label{Xdex}
\\  
&\hbox{where\ } b=\sum_{i=1}^n l_i \om_i\in P,\ j \in 
\frac{1}{ m}\Z,\ \hw\in \hW.
\notag \end{align}
We set $(\tilde{b},\tilde{c})=(b,c)$ ignoring the affine extensions.

Later $Y_{\tb}=Y_b q^{-j}$ will be needed. Note the
negative sign of $j$. 

\begin{definition}
The  double  affine Hecke algebra $\HH\ $
is generated over $ \Q_{ q,t}$ by 
the elements $\{ T_i,\ 0\le i\le n\}$, 
pairwise commutative $\{X_b, \ b\in P\}$ satisfying 
(\ref{Xdex}),
and the group $\Pi,$ where the following relations are imposed:

(o)\ \  $ (T_i-t_i^{1/2})(T_i+t_i^{-1/2})\ =\ 
0,\ 0\ \le\ i\ \le\ n$;

(i)\ \ \ $ T_iT_jT_i...\ =\ T_jT_iT_j...,\ m_{ij}$ 
factors on each side;

(ii)\ \   $ \pi_rT_i\pi_r^{-1}\ =\ T_j \iif 
\pi_r(\al_i)=\al_j$; 

(iii)\  $T_iX_b T_i\ =\ X_b X_{\al_i}^{-1} \iif 
(b,\al^\vee_i)=1,\
0 \le i\le  n$;

(iv)\ $T_iX_b\ =\ X_b T_i$ if $(b,\al^\vee_i)=0 
\for 0 \le i\le  n$;

(v)\ \ $\pi_rX_b \pi_r^{-1}\ =\ X_{\pi_r(b)}\ ,\  r\in O'$.
\label{double}
\end{definition}
\qed

Given $\tw \in \tW, r\in O,\ $ the product
\begin{align}
&T_{\pi_r\tw}\equal \pi_r\prod_{k=1}^l T_{i_k},\where 
\tw=\prod_{k=1}^l s_{i_k},
l=l(\tw),
\label{Twx}
\end{align}
does not depend on the choice of the reduced decomposition
(because $\{T\}$ satisfy the same ``braid'' relations 
as $\{s\}$ do).
Moreover,
\begin{align}
&T_{\hv}T_{\hw}\ =\ T_{\hv\hw}\  \hbox{ whenever}\ 
 l(\hv\hw)=l(\hv)+l(\hw) \for
\hv,\hw \in \hW. \label{TTx}
\end{align}
In particular, we arrive at the pairwise 
commutative elements 
\begin{align}
& Y_{b}\ =\  \prod_{i=1}^nY_i^{l_i} \iif  
b=\sum_{i=1}^n l_i\om_i\in P,\where  
 Y_i\equal T_{\om_i},
\label{Ybx}
\end{align}
satisfying the relations
\begin{align}
&T^{-1}_iY_b T^{-1}_i\ =\ Y_b Y_{\al_i}^{-1} \iif 
(b,\al^\vee_i)=1,
\notag\\ 
& T_iY_b\ =\ Y_b T_i \iif (b,\al^\vee_i)=0,
 \ 1 \le i\le  n.
\end{align}

\smallskip 
The Demazure-Lusztig
operators
are defined as follows: 
\begin{align}
&T_i\  = \  t_i ^{1/2} s_i\ +\ 
(t_i^{1/2}-t_i^{-1/2})(X_{\al_i}-1)^{-1}(s_i-1),
\ 0\le i\le n,
\label{Demazx}
\end{align}
and obviously preserve $\Q[q,t^{\pm 1/2}][X]$.
We note that only the formula for $T_0$ involves $q$: 
\begin{align}
&T_0\  =  t_0^{1/2}s_0\ +\ (t_0^{1/2}-t_0^{-1/2})
( q X_{\vth}^{-1} -1)^{-1}(s_0-1),\notag\\ 
&\where
s_0(X_b)\ =\ X_bX_{\vth}^{-(b,\vth)}
 q^{(b,\vth)},\ 
\al_0=[-\vth,1].
\end{align}

The map sending $ T_j$ to the formula in
(\ref{Demazx}), and $\ X_b \mapsto X_b$ 
(see (\ref{Xdex})),
$\pi_r\mapsto \pi_r$ induces a 
$ \Q_{ q,t}$-linear 
homomorphism from $\HH\, $ to the algebra of linear endomorphisms 
of $\Q_{ q,t}[X]$.
This $\HH\,$-module, which will be called the
{\dfont polynomial representation}, 
is faithful 
and remains faithful when   $q,t$ take  
any nonzero complex values assuming that
$q$ is not a root of unity. 

The images of the $Y_b$ are called the 
{\dfont difference Dunkl operators}. To be more exact,
they must be called difference-trigonometric Dunkl operators,
because there are also 
difference-rational Dunkl operators.

The polynomial representation
is the $\HH\,$-module induced from the one dimensional
representation $T_i\mapsto t_i^{1/2},\,$ $Y_i\mapsto Y_i^{1/2}$
of the affine Hecke subalgebra $\h_Y=\lan T,Y\ran.$
Here the PBW-Theorem is used:
for arbitrary nonzero $q,t,$ any element $H \in \HH\ $  
has a unique decomposition in the form
\begin{align}
&H =\sum_{w\in W }\,  g_{w}\, f_w\, T_w,\ 
g_{w} \in \Q_{q,t}[X],\ f_{w} \in \Q_{q,t}[Y].  
\label{hatdecx}
\end{align}

\smallskip
The definition of DAHA and the polynomial representation
are compatible with the  
{\dfont intermediate subalgebras} $\HH^\flat\subset\HH\ $
with $P$ replaced by any lattice $B\ni b$ between $Q$ and $P.$
Respectively, $\Pi$ is 
changed to the image $\Pi^\flat$ of $B/Q$ in $\Pi.$ 
From now on, we take  $X_a,Y_b$ with
the indices $a,b\in B.$ We will continue using
the notation $\v$ for the $B$-polynomial representation:
$$
\v\ =\ \Q_{q,t}[X_b]\ =\ \Q_{q,t}[X_b, b\in B].
$$

We also set $\hW^{\flat}=B\cdot W\subset \hW,\ $
and replace $m$ by the least $\tilde{m}\in \N$ such that 
$\tilde{m}(B,B)\subset \Z$ in the definition of the 
$\Q_{q,t}.$ 

\medskip
{\bf Automorphisms.}
The following {\dfont duality anti-involution} is of key importance
for the various duality statements: 
\begin{align}
&\phi:\ 
X_b\mapsto Y_b^{-1},\ 
T_i\mapsto T_i\ (1\le i\le n),\label{starphi}.
\end{align}
It preserves $q,t_\nu$ and their fractional
powers. 
 
We will also need the 
automorphisms of
$\HH^{\flat} $(see \cite{C4},\cite{C12}):
\begin{align}
\tau_+:\  &X_b \mapsto X_b, \ Y_r \mapsto 
X_rY_r q^{-\frac{(\om_r,\om_r)}{2}},\
\pi_r \mapsto q^{-(\om_r,\om_r)}X_r\pi_r,
\notag\\ 
\tau_+:\ &Y_\vth \mapsto q^{-1}\,X_\vth T_0^{-1} 
T_{s_\vth},\, T_0\mapsto q^{-1}\,X_\vth T_0^{-1}, 
\and 
\label{taux}\\
\tau_-\ & \equal  \phi\tau_+\phi,\  \ 
\si\equal \tau_+\tau_-^{-1}\tau_+\ =\
\tau_-^{-1}\tau_+\tau_-^{-1}, 
\label{tauminax}
\end{align}
where $r\in O'.$ 
They fix
$T_i\,(i\ge 1),\ t_\nu,\ q$
and fractional powers of $t_\nu,\ q.$
Note that $\tau_-=\si\tau_+\si^{-1}.$

In the definition of $\tau_\pm$ and $\si,$ 
we need to add $q^{\pm 1/(2m)}$ to
$\Q_{q,t}.$

The automorphism $\tau_-$ acts trivially
on $\{T_i(i\ge 0),\,\pi_r\,,Y_b\}.$ Hence 
it naturally acts in the polynomial representation
$\v.$
The automorphism $\tau_+$ and therefore $\si$ do not act in $\v.$
The automorphism $\si$ sends $X_b$ to $Y_b^{-1}$ and 
is associated with the Fourier transform in the DAHA theory.

Actually, all these
automorphisms act in the central extension of the
{\em elliptic
braid group} defined by the relations of $\HH\ ,$
where the quadratic relation is dropped. 
The central extension is by the fractional powers of $q.$ 

\smallskip
The elements $\tau_\pm$ generate the projective $PSL(2,\Z),$
which is isomorphic to the braid group $B_3$ due to Steinberg.
\medskip

\section {Macdonald polynomials}
This definition is due to Macdonald (for
$k_{\sht}=k_{\lng}\in \Z_+ $),
who extended in \cite{M4}
Opdam's nonsymmetric polynomials introduced
in the differential case in \cite{O2}
(Opdam mentions Heckman's contribution in \cite{O2}). 
The general case was considered in \cite{C4}. 

We continue using the same notation $X,Y,T$
for these operators acting in the polynomial
representation. The parameters $q,t$ are generic
in the following definition. 

\begin{definition}\label{YONE}
The {\dfont nonsymmetric Macdonald polynomials} 
$\{E_b, b\in P\}$
are unique (up to proportionality) eigenfunctions of
the operators $\{L_f\equal f(Y_1,\cdots, Y_n), 
f\in \Q[X]\}$
acting in $\Q_{q,t}[X]:$
\begin{align}
&L_{f}(E_b)\ =\ f(q^{-b_\#})E_b, \where
b_\#\equal b- u_b^{-1}(\rho_k),
\label{Yone}
\\ 
& X_a(q^{b})\ =\
q^{(a,b)}\for a,b\in P,\ 
u_b=\pi_b^{-1}b, \label{xaonb}
\end{align}
where $u_b$ is from Proposition \ref{PIOM}. 

They satisfy
\begin{align}
&E_b-X_b\ \in\ \oplus_{c\succ b}\Q(q,t) X_c,\
\langle E_b, X_{c}\rangle_\circ = 0 \for P \ni c\succ b,
\label{macd}
\end{align}
where we set $c\succ b$ if 
$$
c_- -b_-\in B\cap Q_+ \hbox {\ \ or\ \ }
c_-=b_-\and c-b\in B\cap Q_+.
$$
\end{definition}

The following {\dfont intertwiners} 
are the key in the theory:
\begin{align}\label{tauintery}
\Psi_i=\tau_+(T_i)+
\frac{t_i^{1/2}-t_i^{-1/2}}{Y_{\al_i}^{-1}-1},\ i\ge 0,\ \   
P_r=\tau_+(\pi_r),\ r\in O',\\
\Psi_{\hw}=P_r\Psi_{i_l}\ldots\Psi_{i_1}\hbox {\ for\  
reduced\ decompositions\ }
\hw=\pi_rs_{i_l}\ldots s_{i_1}.\notag
\end{align} 
Note the formulas  
$$
\tau_+(T_0)=X_0^{-1}T_0^{-1},\ X_0=qX_{\vth}^{-1},\  
\tau_+(\pi_r)=q^{-(\om_r,\om_r)/2}X_r\pi_r.
$$
The products $\Psi_{\hw}$ do not depend on the
choice of the reduced decomposition,
intertwine $Y_b,$ and transform the $E$-polynomials
correspondingly. Namely, for $\hw\in\hW,$
\begin{align}\label{interhwy}
&\Psi_{\hw}Y_b=Y_{\hw(b)}\Psi_{\hw}, \hbox{\ where\ }
Y_{[b,j]}\equal Y_b q^{-j},\\
&E_{b}=\hbox{Const\,}\Psi_{\hw} (E_c)
\for \hbox{\ Const\ }\neq 0,\ b= \hw\llb c\rrb,\notag
\end{align} 
provided that
$\pi_b=\hw\pi_c$ and $l(\pi_b)=l(\hw)+l(\pi_c).$

Here we use 
the {\dfont affine action} of $\hW$ on 
$z \in \R^n:$
\begin{align}
& (wb)\llb z \rrb \ =\ w(b+z),\ w\in W, b\in P,\notag\\ 
& s_{\tal}\llb z\rrb\ =\ z - ((z,\al)+j)\al,
\ \tal=[\al,\nu_\al j]\in \tR.
\label{afaction}
\end{align}
\smallskip

The definition of the $E$-polynomials and the
action of the intertwiners are compatible
with the transfer to the intermediate subalgebras
$\HH^{\flat}.$ 
Recall that the $B$-polynomial representation is
$$
\v\ =\ \Q_{q,t}[X_b]\equal\Q_{q,t}[X_b,b\in B].
$$
We note that the $\Psi$-intertwiners were introduced by
Knop and Sahi in the case of $GL_n.$ 
\medskip

The coefficients of the Macdonald
polynomials are rational functions in terms of $q_\nu,t_\nu.$
The following evaluation formula holds:
\begin{align}
E_{b}(q^{-\rho_k}) \ =&\ q^{(\rho_k,b_-)} 
\prod_{[\al,j]\in \la'(\pi_b)}
\Bigl(
\frac{ 
1- q_\al^{j}t_\al X_\al(q^{\rho_k})
 }{
1- q_\al^{j}X_\al(q^{\rho_k})
}
\Bigr),
\label{ebebs}\\ 
\la'(\pi_b)\ =&\ 
\{[\al,j]\ |\  [-\al,\nu_\al j]\in \la(\pi_b)\}.
\label{jbseto}
\end{align}
Explicitly,
(see (\ref{lambpi})),
\begin{align}
\la'(\pi_b)\ =& \{[\al,j]\, \mid\,\al\in R_+,
\label{jbset}  
\\  
&-( b_-, \al^\vee )>j> 0 
\iif  u_b^{-1}(\al)\in R_-,\notag\\ 
&-( b_-, \al^\vee )\ge j > 0 \iif   
u_b^{-1}(\al)\in R_+ \}. 
\notag \end{align}
Formula (\ref{ebebs}) is the Macdonald 
{\dfont evaluation conjecture}
in the nonsymmetric variant from \cite{C4}. 

Note that one has to consider only long $\al$ (resp., short)
if $k_{\sht}=0$ (resp., 
$k_{\lng}=0$) in the $\la'$-set.  

We have the following {\dfont duality formula}
for $b,c\in P\, :$
\begin{align}
&E_b(q^{c_{\#}})E_c(q^{-\rho_k})\ =\ E_c(q^{b_{\#}})
E_b(q^{-\rho_k}),\ 
b_\# = b- u_b^{-1}(\rho_k).
\label{ebdual}  
\end{align} 
See \cite{C4}. The proof is based on 
the anti-involution $\phi$ from
(\ref{starphi}).
\medskip

{\bf The action of $\tau_-.$}
The authomorphism
$\tau_+$ is a {\em formal} conjugation
by the Gaussian $\ga(q^z)=q^{(z,z)/2}$ where
we set $X_b(q^z)=q^{(b,z)}.$ We treat $\ga$ as
an element in a completion of the polynomial representation
with the extended action of $\HH.$ Actually, only 
the $W$-invariance of $\ga$ and the relations
$$
\om_j(\ga)=q^{(\om_i,\om_i)/2}X_i^{-1}\ga \for j=1,\ldots,n
$$ 
are needed here.
For instrance, one can (formally) take 
\begin{align}\label{gauserx}
&\ga_x\equal \sum_{b\in B} q^{-(b,b)/2}X_b.
\end{align}

Applying $\si$ and using that $\tau_-=\si\tau_+^{-1}\si^{-1},$ we
obtain that the automorphism $\tau_-$
in $\v$ is proportional to the {\em multiplication}
by 
\begin{align}\label{gausery}
&\ga_y\equal \sum_{b\in B} q^{(b,b)/2}Y_b
\end{align}
provided that $|q|<1.$
We use that $\v$ is a union of finite dimensional spaces
preserved by the $Y$-operators.
This observation
is convenient, although not absolutely necessary, to check
the following proposition.

\begin{proposition}\label{TAUMINE}
i) For generic $q,t$ or for any $q,t$ provided that the polynomial
$E_b$ for $b\in B$ is well-defined, 
\begin{align}\label{taumineb}
\tau_-(E_b)=q^{-\frac{(b_-\,,\,b_-)}{2}-(b_-\,,\,\rho_k)}\,E_b
\for P_-\ni b_-\in W(b). 
\end{align}

ii) For arbitrary $q,t,$ 
\begin{align}
&\tau_-(T_i)=T_i,\ \tau_-(\Psi_i)=\Psi_i \for i>0,\\
&\tau(\tau_+(\pi_r)=q^{(\om_r,\om_r)/2}Y_r\tau_+(\pi_r),\ 
%=q^{-(\om_r,\om_r)/2}\tau_+(\pi_r)Y_{r^*}
\tau_-(\tau_+(T_0))=\tau_+(T_0)^{-1}Y_0,
\notag\\
&\tau_-(\Psi_0)=\Psi_0Y_0=Y_0^{-1}\Psi_0,\ 
Y_0=q^{-1}Y_{\vth}^{-1}.\notag
\end{align}

iii) If $q$ is not a root of unity and $t_\nu$ are arbitrary, then
$\tau_-$ preserves any $Y$-submodule of $\v.$ 
\end{proposition}

{\em Proof.} The first two claims are straightforward.
As for (iii), since $q$ is generic
one can assume that $0<q<1$ and define
$\widetilde{\tau_-}$ as the operator of multiplication
by  $C^{-1}\ga_y$ using (\ref{gausery}) and taking
$$
C=\sum_{b\in B} q^{(b,b)/2}Y_b(1)=
\sum_{b\in B} q^{(b,b)/2}q^{(b,\rho_k)}.
$$
Then $\widetilde{\tau_-}$ coincides with $\tau_-$ 
for generic $k,$ when all $E$-polynomials exist and the
$X$-spectrum of $\v$ is simple, due to (i).
This gives the coincidence
for any $k.$
\sq

\medskip
\section{The radical}
\setcounter{equation}{0}
Following \cite{C3,C4},
we set
\begin{align}
 &\{f,g\}\ =\ \{L_{\imath(f)}(g(X))\}\ =\  
 \{L_{\imath(f)}(g(X))\}(q^{-\rho_k}) \for
f,g\in \v,
\label{symfg}
\\ 
&\imath(X_b)\ =\ X_{-b}\ =\ X_b^{-1},\ 
\imath(z)\ =\  z \for 
z\in \Q_{q,t}\, ,
\notag \end{align}
where $L_f$ is from Definition \ref{YONE}. 
It induces the 
$\Q_{q,t}$-linear anti-\-involution $\phi$ of
$\HH^\flat$ from (\ref{starphi}). 

\begin{lemma}\label{RADI}
For arbitrary nonzero $q,t_{\sht},t_{\lng},$
\begin{align}
&\{f,g\}=\{g,f\} \and
 \{H(f),g\}=\{f,H^\phi\,
(g)\},\ H\in \HH^\flat.
\label{syminv}
\end{align}
The quotient $\v'$ of $\v$ by the radical 
$Rad\equal$Rad$\{\, ,\,\}$ 
of the pairing $\{\ ,\ \}$ is an 
$\HH^\flat$-module such that

a) all $Y$-eigenspaces of $\v'$ are zero or 
one-dimensional,

b) $E(q^{-\rho_k})\neq 0$ if the image $E'$ of 
$E$ in  $\v'$
is a nonzero $Y$-eigenvector.

The radical $Rad$ is the greatest $\HH^\flat$-submodule
in the kernel of the map 
$f\mapsto \{f,1\}=f(q^{-\rho_k}).$ 
\end{lemma}
{\em Proof.} Formulas (\ref{syminv}) are from
Theorem 2.2 of \cite{C4}. Concerning the rest,
let us recall the argument from \cite{C12}.
Since  Rad$\{\, ,\,\}$  is a submodule, the 
form  $\{\ ,\ \}$
is well defined and nondegenerate on $\v'.$
For any pullback $E\in \v$ of $E'\in \v',$
$E(q^{-\rho_k})=\{E,1\}=$ $\{E',1'\}.$
If $E'$ is a $Y$-eigenvector in $\v'$ and 
$E(q^{-\rho_k})$ vanishes , then
$$
\{Q_{q,t}[Y_b](E'),\h_Y^\flat(1')\}=0= 
\{E',\v\cdot\h_Y^\flat(1')\}.
$$
Therefore $ \{E',\v'\}=0,$ which is impossible.
\sq

In the following lemma, $q$ is generic,
but $t_\nu$ are not supposed generic.
The Macdonald polynomials $E_{b}$ always exist for 
$b=b^o,$  satisfying the conditions
\begin{align}\label{minimeb}
&q^{-a_\#}\neq q^{-b^o_\#} \hbox{\ \ for\ all\ \ } a\succ b^o.
\end{align} 
We call such $b^o$ {\dfont primary.}
Sufficiently big $b$ are primary.

\begin{lemma}\label{EINRAD}
i) A $Y$-eigenvector $E\in \v$ belongs to $Rad$ if and only
if $E(q^{-\rho_k})=0.$ 
The equality $E(q^{-\rho_k})=0$ automatically results in
the equalities 
\begin{align}
E(q^{-b^o_\#})=0 \hbox {\ \ for\ all\ \ }
b^o \in B^\star\equal \{b^o\in B\, 
\mid \,  E_{b^o}(q^{-\rho_k})\neq 0\}.
\end{align}

ii) Let us assume that the radical is nonzero.
Then for any constant $C>0\, (1\le i\le n),$
there exists primary $b^o$ such that
$(\al_i,b^o)>C$ and $E_{b^o}(q^{-\rho_k})=0,$
i.e., $E_{b^o}\in Rad.$
\end{lemma}

{\em Proof.} The first claim follows from
Lemma \ref{RADI}. 
If $E\in Rad$ and there is no such $b^o$ for certain $C,$ then
the number of common zeros
of the translations $c(E)$ of $E$ for {\em any} number
of $c\in B$ is infinite, which is impossible because the
degree of $E$ is finite.
\sq

We come to the following theorem generalizing the
description of singular $k$ from \cite{O6}.

\begin{theorem}\label{THMRADI}
Assuming that $q$ is generic,
the radical vanishes if and only if 
$E_{b^o}(q^{-\rho_k})\neq 0$ for all sufficiently big 
primary $b^o,$ i.e., if the product in
the right-hand side of 
(\ref{ebebs}) is nonzero for all $b\in B$
with sufficiently big $(b,\al_i)$ for $i>0.$\sq
\end{theorem}
\medskip

We can define {\dfont quasi-perfect representations} as
$\HH^\flat$-modules which have a nondegenerate
form $\{\ ,\ \}$ satisfying
(\ref{syminv}).
Then the greatest quasi-perfect quotient of the
polynomial represntation is  $\v/Rad.$
Indeed, any quasi-perfect quotient $V$ of $\v$ supplies it with
a form $\{f,g\}_V$ $=\{f',g'\}$ for the images $f',g'$ of
$f,g$ in $V.$ Then a proper linear combination $\{\,,\,\}_o$
of $\{\,,\,\}$ and $\{\,,\,\}_V$ will satisfy $\{1,1\}_o=0,$
which immediately makes it zero identically.

\setcounter{equation}{0}
\section{The irreducibility}
In this section $q,t$ are arbitrary nonzero, including
roots of unity.

\begin{theorem}\label{THMIRRE}
i) If the quotient $\v'$ of the polynomial representation
$\v$ by the radical Rad$\{\, ,\,\}$ is finite dimensional
and $\tau_-$-invariant, then it is an irreducible $\HH^\flat$-module.
The radical is always $\tau_-$-invariant if $q$ is not a root of unity.

ii) At roots of unity, the $\tau_-$-invariance holds
when the radical is $\HH^\flat$-generated by
linear combinations $\sum c_bE_b$ (provided
that $E_b$ exist) over $b$ with coinciding
$q^{-(b_-,b_-)/2-(b_-,\rho_k)}$ from (\ref{taumineb}). 
\end{theorem}

{\em Proof.} 
Using $\phi\tau_-\phi=\tau_+,$
the relation
$$
\{\tau_+f,g\}\ =\ \{f,\tau_-g\} \for f,g\in \v'
$$
defines the action of $\tau_+$ in $\v'$ and therefore 
the action of $\si$ there
satisfying 
$$\tau_+\tau_-^{-1}\tau_+\ =\ \si\  =\ 
\tau_-^{-1}\tau_+\tau_-^{-1}.
$$
The pairing $\{f,g\}_\si\equal\{\si f,g\}=\{f,\si^{-1}g\}$
corresponds to the {\em anti-involution}
$\heartsuit=\si\cdot\phi=\phi\cdot\si^{-1}$ 
of $\HH^\flat,$
sending
\begin{align}\label{heartsuu}
&\heartsuit:\ 
 T_i\mapsto T_i,\, \pi_r\mapsto\pi_r, \,
Y_b\mapsto Y_b,
\, X_b\mapsto T_{w_0}^{-1}X_{\varsigma(b)}T_{w_0}
\end{align}
for $0\le i\le n,\ b\in B.$
  
It holds in either direction, from $f$ to $g$ and
the other way round, but the form $\{f,g\}_\si,$
generally speaking, could be non-symmetric.
Actually it is symmetric, but we do not need it
for the proof.

Using this {\em non-degenerate}  pairing, we proceed
as follows. Any proper $\HH^\flat$-submodule
$\v''$ of $\v'$ contains at least one $Y$-eigenvector $e''$, so we
can assume that $\v''=\HH^\flat e''.$ The corresponding
eigenvalue cannot coincide with that of $1$ thanks to the
previous lemma. Therefore $\{1',\v''\}_\si=0$
for the image $1'$ of $1$ in $\v',$ and the orthogonal
complement of $\v''$ in $\v'$ 
is a {\em proper} $\HH^\flat$-submodule of $\v'$
containing $1',$  which is impossible.

Using Proposition \ref{TAUMINE}, we obtain (ii).
\sq
\medskip

Let us check that the pairing $\{f,g\}_\si$ is symmetric.
First of all,
\begin{align*}
&\{\tau_+(1'),1'\}=\{1',\tau_-(1')\}=\{1',1'\}=1 \Rightarrow\\
&\{1',1'\}_\si=\{\si(1'),1'\}=\{\tau_+(1'),\tau_-(1')\}=
\{\tau_+(1'),1'\}=1.
\end{align*}
Then,
$\{1',f\}_\si-\{f,1'\}_\si=\{(\si-\si^{-1})(1'),f\}=
\{(1-\si^{-2})(1'),f\}_\si.$
However $\si^{-2}$ coincides with $T_{w_o}$
up to proportionality 
in {\em irreducible} $\HH\,$-modules
where $\si$ acts (see \cite{C12}). 
Thus $(1-\si^{-2})(1')$ is proportional
to $1'$ and  must be zero in $\v'$ due to the calculation above.
We obtain that $1'$ is in the radical of the pairing
$\{f,g\}_\si-\{g,f\}_\si,$ which makes this
difference identically zero since $1'$ is a
generator.  
\medskip

The quotient $\v'$ {\em is not} $\tau_-$-invariant
if $q$ is a root of unity and $k$ are {\em generic}.
In this case (see \cite{C4,C12}), all $E_b$
and $\e_b=E_b/E_b(q^{-\rho_k})$ are well defined.
The radical is linearly
generated by the differences $\e_b-\e_c$ when 
$$
u_b=u_c,\ b_-=c_- \mod NA\cap B \for (A,B)=\Z,\ q^N=1.
$$ 
The polynomials $E_b$ and $\e_b$ are $\tau_-$-eigenvectors. Their
eigenvalues are $q^{-(b_-,b_-)/2-(b_-,\rho_k)}.$
Therefore $\tau_-$ does not preserve the radical.
\medskip

{\bf An example of reducible $\v'.$}
For the root system $B_n (n>2),$ let 
\begin{align*}
&n\ge l>n/2+1,\ r=2(l-1),\  k_{\lng}=-\frac{s}{r},\
l,s\in \N,\, (s,r)=1.
\end{align*}
We will assume that $k_{\sht}$ is generic.

Then Theorem \ref{THMRADI}
readily gives that the radical
is zero. Indeed, the numerator of the formula from
(\ref{ebebs}) is nonzero for {\em all} $b$ because
\begin{align}\label{notev}
&q_\al^jt_\al X_\al(q^{\rho_k})=
q_\al^{j+k_\al+(\al^\vee,\rho_k)}\neq 1
\hbox{\ \ for\ any\ \ }\al\in R_+,j>0,
\end{align}
and the denominator is nonzero because
\begin{align}\label{nooneev}
&q_\al^jX_\al(q^{\rho_k})=
q_\al^{j+(\al^\vee,\rho_k)}\neq 1
\hbox{\ \ for\ any\ \ }\al\in R_+,j>0.
\end{align}
We use that $(\al^\vee,\rho_k)$ involves
$k_{\sht}$ unless $\al$ belongs to the root subsystem
$A_{n-1}$ formed by $\ep_l-\ep_m$ in the notation
of \cite{Bo}.

Thus all Macdonald
polynomials $E_b$ are well-defined
and the $Y$-action in $\v$ is semisimple.
The semisimplicity results from (\ref{nooneev}).
\smallskip

The following relation holds:
\begin{align}\label{yesnegev}
&q_\al^jt_\al^{-1}X_\al(q^{\rho_k})=
q_\al^{j-k_\al+(\al^\vee,\rho_k)}= 1
\for \al=\ep_{l},\, j=2(l-1)s
\end{align}
in the notation from \cite{Bo}.
Indeed, $(\al^\vee,\rho_k)=k_{\sht}+2(l-1)k_{\lng}.$
Let 
$$
\tal^\bullet=[-\al,\nu_{\al}j]=[-\ep_l,2(l-1)s].
$$ 
Here $\al$ is short, so $\nu_{\al}=1.$

\begin{proposition}
The polynomial representation
has a proper submodule $\v^\bullet$ 
which is the linear span of
$E_b$ for $b$ such that $\la(\pi_b)$ contains 
$\tal^\bullet.$ The quotient $\v/\v^\bullet$ is irreducible.  
\end{proposition}
{\em Proof.} This statement follows from the
Main Theorem of \cite{C12}. It is easy to check
it directly using the intertwiners from
(\ref{interhwy}). Indeed, given $b,$ the linear span
$\sum_{\hw}\Psi_{\hw}(E_b)$ is an $\HH^\flat$-submodule
of $\v$ when {\em all} $\hw\in \hW$ are taken, not only
the ones satisfying $l(\hw\pi_b)=l(\hw)+l(\pi_b).$
If $\pi_b$ contains $\tal^\bullet$ 
but $\hw\pi_c$ does not, then  $\Psi_{\hw}(E_b)=0$
because the product $\Psi_{\hw}\Psi_{\pi_b}(1)$ 
can be transformed using the homogeneous Coxeter
relations to get the combination 
$$
\cdots(\tau_+(T_i)-t_i^{1/2})(\tau_+(T_i)+t_i^{-1/2})\cdots(1)
$$
somewhere. This combination is identically zero. 
\sq 
\medskip

\section {A non-semisimple example}
\setcounter{equation}{0}
Let us consider the case of $A_1$ assuming that
$q^{1/2}$ is a primitive $2N$-th root of unity.
We set $t=q^k,$ 
$$
B=P=\Z,\ Q=2\Z,\ X=X_{\om_1},\ Y=Y_{\om_1},\ T=T_1.
$$
Thus the $E$-polynomials will be numbered by integers, and
$Y(E_m)=q^{\la_m}E_m$ for 
$$
\la_m=-m_\#,\ m_\#\equal ({m+\sgn(m)k})/2,\  \sgn(0)=-1,
$$ 
provided that $E_m$ exists. The $\la_m$ are called 
weights of $E_m.$

Note that $\pi=sp$ in the polynomial representation
$\v=\Q_{q,t}[X,X^{-1}]$ for
$\ s(f(X))=f(X^{-1}),\ $ $p(f(X))\equal f(q^{1/2}X).$
The definition ring is 
$\Q_{q,t}=\Q[q^{\pm 1/4},t^{\pm 1/2}],$ where $q^{1/4}$
is use to introduce of $\tau_{\pm}.$ Otherwise $q^{1/2}$
is sufficient.
 
We will need the following lemma, which is 
similar to the considerations from \cite{CO}.

Let $\widehat{V}_0=\Q_{q,t},$ 
$\widehat{V}_1=\Q_{q,t} X,\ldots,$  
$$
\widehat{V}_{-m}=B_m\widehat{V}_m,\
\widehat{V}_{m+1}=A_{-m}\widehat{V}_{-m},\ldots,
$$
where $m>0$, $A_{-m}=q^{m/2}X\pi,$
$B_m$ is the restriction
of the intertwiner 
$t^{1/2}(T+\frac{t^{1/2}-t^{-1/2}}{Y^{-2}-1})$ 
to $\widehat{V}_m$ provided that  $q^{2\la_m}\neq 1$
for $\la_m=-m/2-k/2.$ If $q^{2\la_m}=1$ and
the denominator of $B_m$ becomes infinity, then we set
$\ B_m=t^{1/2}T,\ $ $\ \widehat{V}_{-m}=\widehat{V}_m+
T\widehat{V}_m.\ $

\begin{lemma}\label{epolyn}
i) The space $\widehat{V}_{\pm m}$ is one-dimensional or
two-dimensional. In the latter case, it is the Jordan $2$-block 
satisfying $(Y-q^{\pm\la_m})^2\widehat{V}_{\pm m}=\{0\}.$ 
If dim$\widehat{V}_{-m}=1$ then dim$\widehat{V}_{m+1}=1$
and the generators are
$$E_{-m}=B_{m-1}\cdots B_1A_0(1),\ E_{m+1}=A_{-m}E_{-m}.$$  
If dim$\widehat{V}_{-m}=2,$ then
dim$\widehat{V}_{m+1}=2$ and the $E$-polynomials
$E_{-m},E_{m+1}$ do not exist, although these
spaces contain the $E$-polynomials of smaller degree.

ii) Let us assume that
either $q^{2\la_m}=t$ or  $q^{2\la_m}=t^{-1}.$ 
Then dim$\widehat{V}_{-m}=1$ and this space is
generated by $E_{-m}.$ If $\widehat{V}_{m}$ is 
one-dimensional then respectively 
$(T+t^{-1/2})E_{-m}=0$  or  
$(T-t^{1/2})E_{-m}=0.$
If dim$\widehat{V}_{m}=2,$ then respectively  
$$ (T+t^{-1/2})E_{-m}\hbox{\ \ or\ \ } (T-t^{1/2})E_{-m}$$ 
is nonzero and proportional to
the (unique) $E$-polynomial which is contained in the space 
$\widehat{V}_{m}.$ \sq
\end{lemma}

We are going to apply the lemma to {\em integral}
$k.$ In the range $0<k<N/2,$ the corresponding 
perfect representation is $Y$-semisimple.
Using the reduction modulo $N$ (see \cite{CO}), it
suffices to consider the interval  $-N/2\le k<0.$

\begin{proposition}
i) For integral $k$ such that $\,-N/2\le k<0,\,$ 
the quotient
$V_{2N+4|k|}\equal\v/Rad$ by the radical of the pairing 
$\{\, ,\, \}$ is an irreducible $\HH\,$-module
of dimension $2N+4|k|.$ 

ii) The polynomials $E_m$ exist and
$E_m(q^{-k/2})\neq 0$ for the sequences: 
\begin{gather*}
m=\ \{0,1,-1,\ldots,-|k|+1,|k|\},\\ 
m=\ \{-2|k|,2|k|+1,\ldots,-N+1,N\}, \\
m=\ \{-N,N+1,\ldots,-N-|k|+1,N+|k|\},
\end{gather*}
respectively with $2|k|,$ $2(N-2|k|),$ and $2|k|$ elements.
They do not exist for $2|k|+2|k|$ indices
\begin{gather*}
m=\ \{-|k|,|k|+1,\ldots,-2|k|+1,2|k|\},\\ 
m=\ \{-N-|k|,N+|k|+1,\ldots,-N-2|k|+1,N+2|k|\}.
\end{gather*}

iii) The $Y$-semisimple component 
of $V_{2N+4|k|}$ of dimension $2N-4|k|$ 
is linearly generated by  $E_m$ for 
$$ \{m=-2|k|,\,2|k|+1,\,-2|k|-1,\,\ldots,\,-N+1,\,N\}. 
$$
The corresponding $Y$-weights are 
$$
\{\la=\frac{|k|}{2},\frac{-|k|-1}{2},\frac{|k|+1}{2},\ldots,
\frac{N-1-|k|}{2},\frac{|k|-N}{2}\}.
$$

iv)
The rest of $V_{2N+4|k|}$
is the direct sum of $4|k|$ Jordan $2$-blocks
of the total dimension $8|k|.$ There are
two series of the corresponding
(multiple) weights $\la:$
$$\{\frac{-|k|}{2},\frac{|k|-1}{2},\ldots,\frac{-1}{2},\frac{0}{2}\},\  
\{\frac{N-|k|}{2},\frac{|k|-N-1}{2},
\ldots,\frac{N-1}{2},\frac{-N}{2}\}.
$$
\end{proposition}
{\em Proof.}
We will use the chain of the spaces of generalized
eigenvectors
$$
\widehat{V}_0=\Q_{q,t},\ 
\widehat{V}_{1}=\Q_{q,t} X,\ \widehat{V}_{-1},\ \ldots,\ 
\widehat{V}_m,\ \ldots 
$$
from Lemma \ref{epolyn}. 
Recall that $m>0.$ The following holds:
\smallskip

\noindent
0) the spaces $\widehat{V}_{\pm m}$ are all one-dimensional
from $0$ to $m=|k|,$ i.e., in the sequence $V_0,\ldots,V_{-|k|+1},V_{|k|};$

\noindent
1) the intertwiner $B_m$ becomes infinity at $m=|k|$ ($B_{|k|}=t^{1/2}T$) 
and dim$\widehat{V}_m =2$ in the range $|k|< m\le 2|k|;$

\noindent
2) the intertwiner $B_m$ kills $1\in \widehat{V}_m$ at $m=2|k|$, and
after this dim$\widehat{V}_m =1$ for $2|k|<m\le N;$ 

\noindent
3) $B_m$ is proportional to $(T+t^{-1/2})$ at $m=N,\ $
$E_{-N}=X^N+X^{-N},\ $ and dim$\widehat{V}_m =1$ as $N< m\le N+|k|;$

\noindent
4) the intertwiner $B_m$ becomes infinity again at $m=N+|k|,$
and afterwards dim$\widehat{V}_m =2$ when $N+|k|< m\le N+2|k|;$

\noindent
5) $B_m$ kills $E_{-N}$ at $m=N+2|k|,$ and $B_m(\widehat{V}_m)$ is 
generated by $E_{-N-2|k|}$ of same $Y$-eigenvalue as $E_{N}.$
\smallskip

Concerning step (5),
the polynomials $E_{-N-2|k|}$ and $E_N$ both exist,
there evaluations are nonzero,
and the difference 
$$
E=E_N/E_N(q^{-k/2})-E_{-N-2|k|}/E_{-N-2|k|}(q^{-k/2})
$$ 
belongs to the radical $Rad,$ i.e., becomes zero in $V_{2N+4|k|}.$  

Note that $(T+t^{1/2})E=0,$ which is important to know to
continue the decomposition of $\v$ further. It follows the
same lines.
 
We see that step (5) is the first step which produces
no new elements in $V_{2N+4|k|}.$ Namely:
$$
B_{N+2|k|}(\widehat{V}_{N+2|k|})\ =\ 
\Q_{q,t}E_N \hbox{\ \ in\ \ } V_{2N+4|k|},
$$
and we can stop here.

The lemma gives
that between (2) and (3), the polynomials $E_m$ exist, 
their images linearly generate the $Y$-semisimple part of $V.$
It is equivalent to the inequalities 
$E_m(q^{-\rho_k})\neq 0$ because
they have different $Y$-eigenvalues.
 
Apart from (2)-(3), there
will be Jordan $2$-blocks with respect to $Y.$
Let us check it.

First, we obtain the  
2-dimensional irreducible representation
of $\h_Y=\langle T,Y,\pi \rangle$ in the corresponding
$\widehat{V}$-space at step (1). Then we apply 
invertible intertwiners to this space
(the weights will go back) and eventually will
obtain the two-dimensional $\widehat{V}$-space 
for the starting weight $\la=-|k|/2.$ 
Note that $E_0=1$ is not from
the $Y$-semisimple component of
$V_{2N+4|k|}.$ It belongs to a Jordan $2$-block.

Second, the intertwiner (2)  makes the last space
one-dimensional and $Y$-semisimple (the 
corresponding eigenvalue is simple in $ V_{2N+4|k|}$).
It will remain one-dimensional
until (3). After step (3), we obtain the Jordan blocks. The steps
(4)-(5) are parallel to  (1)-(2). 
\sq
\medskip

The above consideration
readily results in the irreducibility of 
the module $V_{2N+4|k|}.$  
Indeed, Lemma \ref{epolyn}, (ii) gives
that if a submodule of $V_{2N+4|k|}$
contains at least one simple $Y$-eigenvector
then it contains the image of $1$ and the whole
space. Step (5) guarantees that it is
always the case, because we can obtain $E_N$
beginning with an arbitrary $Y$-eigenvector.
\smallskip

The irreducibility and the existence of
the projective $PSL(2,\Z)$-action in $V_{2N+4|k|}$
also follow from
Theorem \ref{THMIRRE}, (ii) because the radical
is generated by 
$E$ which is a linear combination 
of the $E$-polynomials with the coinciding
$\tau_-$-eigenvalues.

\medskip
\bibliographystyle{unsrt}

\begin{thebibliography} {ABCD}
%\vfil
%\medskip

\bibitem [B] {Bo}
{N.~Bourbaki},
{\em Groupes et alg\`ebres de Lie}, Ch. { 4--6},
Hermann, Paris (1969)

\bibitem [C1] {C3}
{I.~Cherednik},
{\em Macdonald's evaluation con\-jec\-tures and
dif\-fe\-rence Fou\-rier tran\-sform},
Inventiones Math. {122} (1995),119--145.

\bibitem [C2] {C4}
\bysame, 
{\em Nonsymmetric Macdonald polynomials },
IMRN { 10} (1995), 483--515.

\bibitem [C3] {C12}
\bysame, 
{\em Double affine Hecke algebras and difference Fourier transforms},
Inventiones Math. { 152} (2003), 213--303.

\bibitem [C4] {C29}
\bysame, 
{\em Diagonal coinvariants and Double Affine Hecke algebras},
IMRN { 16}  (2004), 769--791.

\bibitem [CO] {CO}
\bysame, and {V.~Ostrik},
{\em From Double Hecke Algebras to
Fourier Transform}, Selecta Math. New ser. { 9} (2003), 1022--182.

\bibitem [DO] {DO}
{C.F.~Dunkl}, and {E.M.~Opdam},
{\em Dunkl operators for complex reflection groups},
Proc. London Math. Soc. (3), { 86} (2003), 70--108. 

\bibitem [DJO] {DJO}
{C.F.~Dunkl}, and {M.~de Jeu}, and {E.M.~Opdam},
{\em Singular polynomials for finite reflection groups},
Trans. Amer. Math. Soc. { 346} (1994), 237--256. 

\bibitem [J] {Je1}
{M.~de Jeu},
{\em The Dunkl operators}, Thesis (1993), 1--92.

\bibitem [M] {M4}
{I.~Macdonald}, 
{\em Affine Hecke algebras and orthogonal polynomials},
S\'e\-mi\-naire Bour\-baki { 47}:797 (1995), 01--18.

\bibitem [O1] {O2}
{E.~Opdam}, 
{\em Harmonic analysis for certain 
representations of graded Hecke algebras}, 
Acta Math. { 175} (1995), 75--121.

\bibitem [O2] {O6}
\bysame,
{\em Dunkl operators, Bessel functions and the
discriminant of a finite Coxeter group},
Composito Mathematica { 85} (1993), 333--373.


\end{thebibliography}

\end{document}